\documentclass{amsart}

\usepackage{amssymb}
\usepackage{amsthm}
\usepackage{amsmath} 
\usepackage{amsbsy}
\usepackage{hyperref}
\usepackage{tikz}
\newtheorem{theorem}{Theorem}

\newtheorem{lemma}{Lemma}

\newcommand{\abs}[1]{\lvert #1 \rvert}

\begin{document}

\title[]{A Dimension-Free Hermite-Hadamard inequality\\ via gradient estimates for the torsion function} \keywords{Hermite-Hadamard, subharmonicity, Brownian motion}
\subjclass[2010]{26B25, 28A75, 31A05, 31B05, 35B50.} 

\thanks{J.L.~is supported in part by the National Science Foundation via grant DMS-1454939. S.S. is supported in part by the NSF (DMS-1763179) and the Alfred P. Sloan Foundation.}

\author[]{Jianfeng Lu}
\address[Jianfeng Lu]{Department of Mathematics, Department of Physics, and Department of Chemistry,
Duke University, Box 90320, Durham NC 27708, USA}
\email{jianfeng@math.duke.edu}

\author[]{Stefan Steinerberger}
\address[Stefan Steinerberger]{Department of Mathematics, Yale University, New Haven, CT 06510, USA}
\email{stefan.steinerberger@yale.edu}

\begin{abstract} 
Let $\Omega \subset \mathbb{R}^n$ be a convex domain and let $f:\Omega \rightarrow \mathbb{R}$ be a subharmonic function, $\Delta f \geq 0$, which satisfies $f \geq 0$ on the boundary $\partial \Omega$. Then
$$ \int_{\Omega}{f ~dx} \leq  |\Omega|^{\frac{1}{n}} \int_{\partial \Omega}{f ~d\sigma}.$$
Our proof is based on a new gradient estimate for the torsion function, $\Delta u = -1$ with Dirichlet boundary conditions, which is of independent interest.
\end{abstract}
\maketitle

\section{Introduction and results}
 \subsection{Introduction.} The Hermite-Hadamard inequality is an elementary observation, attributed to both Hadamard \cite{hada} and Hermite \cite{hermite}, for convex functions on the real line $f:[a,b] \rightarrow \mathbb{R}$ stating that
$$ \frac{1}{b-a} \int_{a}^{b}{f(x) dx} \leq \frac{f(a) + f(b)}{2}.$$
Many generalizations of all sorts are known, we refer to the monograph of Dragomir \& Pearce \cite{drago}. However, there is relatively little work outside of the one-dimensional setting, we refer to \cite{cal1, cal2, nicu2, nicu, choquet, stein}. It is known \cite{stein} that if $\Omega \subset \mathbb{R}^n$ is a convex domain and if $f:\Omega \rightarrow \mathbb{R}$ is a convex function such that $f \big|_{\partial \Omega} \geq 0$, then 
\begin{equation}
 \frac{1}{|\Omega|} \int_{\Omega}{f ~d x} \leq \frac{c_n}{|\partial \Omega|} \int_{\partial \Omega}{f ~d \sigma}
 \end{equation}
for some universal constant $c_n > 0$ depending only on the dimension. It is known that the optimal constant satisfies $1 < c_n \leq 2n^{n+1}$ but it is not clear what type of growth to expect (or whether there is any growth at all). The isoperimetric inequality immediately implies that, for some other universal constants $c_n^*$,
\begin{equation}
 \int_{\Omega}{f ~d x} \leq  c_n^* |\Omega|^{1/n}\int_{\partial \Omega}{f ~d \sigma},
 \end{equation}
 where $c_n^*$ satisfies $c_n^* \lesssim c_n n^{-1/2}$.
The dependence on $\abs{\Omega}^{1/n}$ on the right hand side is consistent with scaling. It is also known \cite{stein} that for $n=2$ it is enough that $\Omega$ is simply connected and that $|\Omega|^{1/2}$ can be replaced by $\mbox{inrad}(\Omega)$ (this is a consequence of the rigidity of $\mathbb{R}^2$, no such results are possible in higher dimensions).

\subsection{A Hermite-Hadamard inequality.} The purpose of our short note is to show that the constants in (2) can be replaced by a universal constant. Moreover, we will prove the inequality for the more general class of subharmonic functions.
\begin{theorem}Let $\Omega \subset \mathbb{R}^n$ be a convex body, $n \geq 2$, and let $f:\Omega \rightarrow \mathbb{R}$ be a subharmonic function, $\Delta f \geq 0$, which satisfies $f \geq 0$ on the boundary $\partial \Omega$. Then,
$$ \int_{\Omega}{f(x)dx} \leq \frac{\sqrt{2}}{\pi} |\Omega|^{1/n} \int_{\partial \Omega}{f d\sigma}.$$
Moreover, the constant $\sqrt{2}/\pi \sim 0.45$ cannot be replaced by $0.22$.
\end{theorem}
This establishes the existence of a universal constant for (2), the analogous question for (1) is still open, though it is certainly conceivable that the optimal constants in (1) do not stay bounded. Likewise one could ask whether one can improve Theorem 1 as the dimension increases (our example showing it cannot be less than 0.22 is in two dimensions), maybe there is actual decay in the optimal constant as the dimension increases? 

\subsection{A gradient estimate.} Theorem 1 follows as a byproduct of a new estimate for an elliptic equation: let $\Omega \subset \mathbb{R}^n$ be a convex domain and consider the equation 
\begin{align*}
- \Delta u &= 1  \qquad \mbox{inside}~\Omega\\
u &= 0\qquad \mbox{on}~\partial \Omega.
\end{align*}
This equation has several different names: the torsion function
\cite{ban} or the Saint Venant problem \cite{carbery} in the theory of elasticity, the landscape
function \cite{svitlana} in the theory of localization of eigenfunctions of elliptic operators and half of the expected lifetime of Brownian
motion before exiting the domain $\Omega$ in probability theory. It is also one of the
elliptic `benchmark' PDEs which are usually the first testcase for new
results \cite{beck, makar, stein1}. We establish a new result for this
equation; it is an interesting question to which extent it can be
generalized to other equations.

\begin{theorem} Let $\Omega \subset \mathbb{R}^n$ be convex, $n \geq 2$, and let $-\Delta u = 1$ with Dirichlet conditions imposed on $\partial \Omega$. Then we have the gradient estimate
  $$\max_{x \in \partial \Omega}{ ~\frac{\partial u}{\partial \nu} } \leq  \frac{\sqrt{2}}{\pi}|\Omega|^{1/n}$$
  where $\nu$ is the interior unit normal on $\partial \Omega$.
\end{theorem}
To the best of our knowledge, no such result is known. Gradient bounds for the torsion function are classical \cite{wh, keady, pa1, pa2}, however, most of them seem to be concerned with the physically relevant case $n=2$. Theorem 2 is sharp up to the value of the constant -- this can be seen in the following example that was suggested to us by Thomas Beck: let 
$$ \Omega = \left\{x \in \mathbb{R}^n: 1 - \frac{x_1^2}{4} - \frac{x_2^2}{2n-2} - \frac{x_3^2}{2n-2} - \dots - \frac{x_n^2}{2n-2} \geq 0 \right\}.$$
The function defining $\Omega$ is already the torsion function on its domain and 
  $$\max_{x \in \partial \Omega}{ ~\frac{\partial u}{\partial \nu} }= 1.$$
 At the same time, we can compute the volume of the ellipsoid and obtain
 $$ | \Omega|^{1/n} = \omega_n^{1/n} 2^{1/n} (2n-2)^{\frac{n-1}{2n}}$$
 which converges to $2\sqrt{\pi e}$ as $n \rightarrow \infty$.  
 
The connection of Theorem 2 to our original problem is via integration by parts which shows that
\begin{align*}
 \int_{\Omega}{f(x) dx} &= \int_{\Omega}{f(x) (-\Delta u(x)) dx} \\
 &= \int_{\Omega}{-\Delta f(x) u(x) dx} + \int_{\partial \Omega}{ \frac{\partial u}{\partial \nu} f d\sigma} \\
 &\leq \int_{\partial \Omega}{ \frac{\partial u}{\partial \nu} f d\sigma} \leq \max_{x \in \partial \Omega}{ ~\frac{\partial u}{\partial \nu} }  \int_{\partial \Omega}{f d\sigma},
\end{align*}
where we used integration by parts ($\nu$ being the \textit{inward} pointing normal vector), subharmonicity $\Delta f \geq 0$ in combination with the maximum principle ($u \geq 0$) and, finally, that $f$ restricted to the boundary is positive $f \big|_{\partial \Omega} \geq 0$ (this idea can be found in a paper of Niculescu and Persson \cite{choquet}).

\section{Some Lemmata}

 We collect several classical results; none of them are new.
\begin{lemma}[Reflection principle; e.g.~\cite{katz}] Let
  $\varepsilon > 0$ and let $T_{\varepsilon}$ be the smallest stopping
  time for a Brownian motion on the real line started at $B(0) = 0$,
  $T_{\varepsilon} = \inf\left\{ t > 0: B(t) = \varepsilon\right\}.$
  Then
$$ \mathbb{P}(T_{\varepsilon} \leq t) = 2 - 2\Phi\left( \frac{\varepsilon}{\sqrt{t}} \right),$$
where $\Phi$ is the c.d.f. of $\mathcal{N}(0,1)$.
\end{lemma}
In particular, this implies that the distribution $\psi$ of the stopping time is given by
$$ \psi(t) = \frac{d}{dt} \mathbb{P}(T_{\varepsilon} \leq t) =  \Phi' \left( \frac{\varepsilon}{\sqrt{t}} \right) \frac{\varepsilon}{t^{3/2}} = \frac{ \varepsilon}{\sqrt{2\pi} t^{3/2}} e^{-\frac{\varepsilon^2}{2t}}.$$
We observe that this distribution decays roughly like $\sim t^{-3/2}$ and does not have a finite mean. However, the expected lifetime of particles hitting the threshold before time $T$ can be computed and is bounded by
$$ \mathbb{E}(T_{\varepsilon} \big| T_{\varepsilon} \leq T) = \int_{0}^{T}{ \psi(t) t dt} \leq \int_0^T { \frac{\varepsilon}{\sqrt{2\pi} t^{1/2}} dt} = \varepsilon \sqrt{\frac{2}{\pi}} \sqrt{T}.$$

\begin{lemma}[see e.g.~\cite{bandle, brock, almut}] Let
  $\Omega \subset \mathbb{R}^n$ be a set with finite volume that
  contains $x_0$ in the interior. Then the expected lifetime of
  Brownian motion started at $x_0$ is maximized when $\Omega$ is a
  ball centered around $x_0$. 
    
\end{lemma}
We want a quantitative estimate with an explicit constant, which is
easy to derive from Lemma 2.  The expected lifetime of a Brownian
motion inside a domain $\Omega$ is given by the torsion function
(note the right hand side is changed to $2$ as the
  infinitesimal generator for standard Brownian motion is
  $\Delta/2$)
\begin{align*}
- \Delta u &= 2  \qquad \mbox{inside}~\Omega\\
u &= 0\qquad \mbox{on}~\partial \Omega
\end{align*}
We will now solve this equation inside a ball. The solution is radial, this motivates the ansatz $u(r) = c(R^2-r^2)$ on a ball of radius $R$ for a constant $c$ to be specified, plugging in shows that $c = n^{-1}.$
Parabolic scaling then shows that Lemma 2 immediately implies the following Lemma 3.

\begin{lemma}[quantitative estimate; see e.g.~\cite{bandle, brock,
    almut}] Let $\Omega \subset \mathbb{R}^n$ be a set with finite
  volume and let $x_0 \in \Omega$. Then the expected lifetime
  $T_{x_0}$ of Brownian motion before hitting $\Omega^c$, started $x_0$, satisfies
$$ \mathbb{E} T_{x_0} \leq \frac{1}{n}\left( \frac{|\Omega|}{\omega_n}\right)^{\frac{2}{n}}.$$
\end{lemma}

\section{Proof}
We prove
  Theorem 2, which implies Theorem 1 by the mechanism outlined
above.  As for Theorem 2, we will work with the explicit interpretation
of the solution of $-\Delta u = 2$ as the expected lifetime of Brownian
motion and derive a bound on that lifetime inside a convex domain
when starting at distance $\varepsilon$ from the boundary. We will
obtain a bound that is linear in $\varepsilon$ and this then implies
the gradient bound by taking a limit.

  \begin{proof}  Our proof of Theorem 2 is based on explicitly
  estimating the expected lifetime for a Brownian motion starting in a
  point $x_0$ located at distance $\varepsilon$ from the boundary
  $\partial \Omega$. We observe that, since $\Omega$ is convex, there
  exists a supporting hyperplane at the point
  $x_1 \in \partial \Omega$ that is distance $\varepsilon$ from $x_0$
  (see Fig.~1). There might be more than one point on the boundary
  that is of minimal distance and in each of these points there might
  be more than one supporting hyperplane -- in case of ambiguity, any
  choice is admissible.  We select a time $T>0$ and study the
  \begin{enumerate}
  \item expected lifetime of particles colliding with the plane within $T$
  units of time,
  \item the likelihood of a particle colliding with the plane within $T$
  units of time, 
  \item and what happens to the surviving articles.
\end{enumerate}

The first step is easy because of domain monotonicity: the likelihood
of a Brownian motion impacting the boundary is made smaller if we
enlarge the domain: in the first step we will thus not use any
information about the convex domain except for the fact that $\Omega$
is in the upper half space defined by the supporting
hyperplane. 
 \begin{center}
\begin{figure}[h!]
\begin{tikzpicture}
\draw [thick] (0,0) -- (10,0);
\draw [thick] (2,1)  to[out=330, in=180] (5,0) to[out=0, in=250] (8.5, 2);
\filldraw (5, 1) circle (0.04cm);
\node at (5.3, 0.8) {$x_0$};
\draw [dashed] (5, 0.8) -- (5, 0);
\node at (4.8, 0.4) {$\varepsilon$};
\filldraw (5, 0) circle (0.04cm);
\node at (5.3, -0.2) {$x_1$};
\node at (2, -0.4) {supporting hyperplane};
\node at (6,2) {$\Omega$};
\end{tikzpicture}
\caption{An outline of the geometry.}
\end{figure}
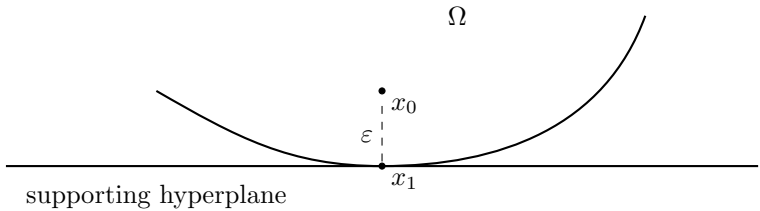
\end{center}
However, that exit time is fairly easy to compute
for the upper half space since Brownian motion is
then independent in all its coordinates and thus the only remaining
question is how long it takes for one-dimensional Brownian motion to
cross a threshold of size $\varepsilon$ (Lemma 1).  We now fix time
$T$.  We have seen in Lemma 1 that the expected lifetime up to time
$T$ satisfies
$$ \int_{0}^{T}{ \psi(t) t\, dt} \leq \varepsilon \sqrt{\frac{2}{\pi}} \sqrt{T}.$$
The probability of surviving (i.e. not colliding with the boundary up to time $T$) is
 $$ \mathbb{P}(T_{\varepsilon} > T) =  2\Phi\left( \frac{\varepsilon}{\sqrt{T}} \right) - 1,$$
 where $\Phi$ is the c.d.f. of the standardized Gaussian $\mathcal{N}(0,1)$.

We use the trivial estimate, valid for all $x>0$,
$$ \Phi\left( x \right)  \leq \frac{1}{2} + \frac{x}{\sqrt{2\pi}}$$
and obtain
 $$ \mathbb{P}(T_{\varepsilon} > T) \leq \sqrt{\frac{2}{\pi}} \frac{\varepsilon}{\sqrt{T}}.$$
 It remains to understand what happens to particles that do not impact within the first $T$ units of time. We do not know where they are after $T$ units of time but we do have a uniform estimate for how much longer we expect them to live (this is Lemma 3).
 Using Markovianity and collecting all these estimates, we thus obtain that
 \begin{align*}
  \mathbb{E}~ \mbox{lifetime} &\leq \varepsilon \sqrt{\frac{2}{\pi}} \sqrt{T} +  \sqrt{\frac{2}{\pi}} \frac{\varepsilon}{\sqrt{T}} \left( T + \frac{1}{n}\left( \frac{|\Omega|}{\omega_n}\right)^{\frac{2}{n}}\right)\\
  &=  \varepsilon \sqrt{\frac{2}{\pi}} \left( 2\sqrt{T} +  \frac{1}{\sqrt{T}} \frac{1}{n}\left( \frac{|\Omega|}{\omega_n}\right)^{\frac{2}{n}}\right).
  \end{align*}
  This quantity is minimized for 
  $$ T= \frac{1}{2n}\left( \frac{|\Omega|}{\omega_n}\right)^{\frac{2}{n}}$$ and results in the bound
 $$  \mathbb{E}~ \mbox{lifetime} \leq   \varepsilon \frac{4}{\sqrt{\pi}} \frac{1}{\sqrt{n}}\left( \frac{|\Omega|}{\omega_n} \right)^{\frac{1}{n}}.$$
 However, this bound is linear in $\varepsilon$ and shows
 $$\max_{x \in \partial \Omega}{ ~\frac{\partial u}{\partial \nu} } \leq \frac{4}{\sqrt{\pi}}\frac{1}{\sqrt{n}}\left( \frac{|\Omega|}{\omega_n} \right)^{\frac{1}{n}}.$$
We assume that $n\geq 2$ in which case we
note the elementary inequality
 $$ \omega_n^{1/n} \sqrt{n} \geq \sqrt{2\pi} \qquad \mbox{for all}~n \geq 2$$
 and thus
  $$\max_{x \in \partial \Omega}{ ~\frac{\partial u}{\partial \nu} } \leq  \frac{\sqrt{8}}{\pi} |\Omega|^{1/n}.$$
  This bound is for the equation $-\Delta u = 2$, the bound for the solution of the equation $-\Delta u = 1$ follows from scaling.
\end{proof}
Moreover, we have that $ \omega_n^{1/n} \sqrt{n} \leq \sqrt{2 \pi e}$ and
$$  \lim_{n \rightarrow \infty} \omega_n^{1/n} \sqrt{n} = \sqrt{2 \pi e}$$
and this shows that the bound $w_{n}^{1/n} \sqrt{n} \geq \sqrt{2\pi}$ is not too lossy.\\

\textbf{Example.}  The example
$$ \Omega = \left\{(x,y): x^2 + y^2 \leq 1 \wedge y \geq 0\right\}$$
and $f(x,y) = 1-y$ show that the constant cannot be replaced by $\sim 0.22$. In general, however, the proof also shows that we get improved estimates in higher dimensions: $\omega_n^{1/n} \sqrt{n}$ converges to (and, in particular, is arbitrarily close to for sufficiently high dimension) $\sqrt{2\pi e}$. We can also consider the radial function
$$ u(r) = \frac{R^2 - r^2}{2n} \qquad \mbox{satisfying} \qquad \Delta u = -1~\mbox{in}~\mathbb{R}^n$$
and consider it on the ball of radius $R$.
The normal derivative pointing in the interior then satisfies
  $$\max_{x \in \partial \Omega}{ ~\frac{\partial u}{\partial \nu} } = \frac{R}{n} \qquad \mbox{while} \qquad |\Omega|^{1/n} =  \omega_n^{1/n}R.$$
This shows that the best constant $c_n$ in the inequality
  $$\max_{x \in \partial \Omega}{ ~\frac{\partial u}{\partial \nu} } \leq c_n |\Omega|^{1/n}$$
  necessarily satisfies
  $$ c_n \geq \frac{1}{n \omega_n^{1/n}} = \frac{1}{\sqrt{n}} \frac{1}{\sqrt{n} \omega_n^{1/n}} \geq \frac{1}{2} \frac{1}{\sqrt{n}}.$$

\end{document}